\documentclass[11pt]{amsart}
\usepackage{amsxtra,amssymb,amsmath,amsthm,amsbsy}
\usepackage[all]{xy}

\usepackage{enumerate}

\numberwithin{equation}{section}

\theoremstyle{plain}
\newtheorem{theorem}{Theorem}[section]

\theoremstyle{definition}
\newtheorem{definition}[theorem]{Definition}
\newtheorem{example}[theorem]{Example}

\def\R{\mathbb R}

\def\then{\Rightarrow}

\def\action{\curvearrowright}
\def\toto{\rightrightarrows}

\def\<{\langle}
\def\>{\rangle}

\newcommand{\Rm}         {(\mathbb{R},\cdot)}

\begin{document}

\title[]
{Lie theory of vector bundles, Poisson geometry and double structures}

\author[]
{Henrique Bursztyn \and Alejandro Cabrera \and Matias del Hoyo}

\address
{IMPA, Estrada Dona Castorina 110, Rio de Janeiro, 22460-320, Brazil
} \email{henrique@impa.br, mdelhoyo@impa.br}

\address{Departamento de Matem\'atica Aplicada - IM, Universidade
Federal do Rio de Janeiro, CEP 21941-909, Rio de Janeiro, Brazil}
\email{acabrera@labma.ufrj.br}

\address{}
\email{}

\date{}


\begin{abstract}
We briefly review our results on the Lie theory underlying vector bundles over Lie groupoids and Lie algebroids, pointing out the role of Poisson geometry in extending these results to double Lie algebroids and LA-groupoids.

\end{abstract}


\maketitle

\vspace{-1cm}

\tableofcontents

\vspace{-1cm}


\section{Lie groupoids and Lie algebroids in mathematics and physics}


Lie groupoids are categorified manifolds which offer a common framework for various classical geometries, such as group actions, fibrations and foliations, codifying internal and external symmetries. Lie groupoids
also play a key role as models for singular spaces, either regarded as generalized atlases (as set forth by Grothendieck), or by providing noncommutative models via their convolution algebras \cite{CW,Connes}. Lie algebroids are their infinitesimal counterparts; they generalize tangent bundles and Lie algebras,
and bear intimate connections with Poisson geometry \cite{CW,CDW}.


The relevance of Lie groupoids and algebroids in mathematical physics is a natural consequence of the deep ties between classical mechanics and Poisson geometry, and between quantum theory and noncommutative geometry; for example,
Lie groupoids and algebroids arise in geometric mechanics \cite{we}, as part of quantization schemes \cite{BatesWein,Hawkings}, and in the theory of sigma models \cite{BKS,catfel}. Lie groupoids form a 2-category which is also central in the study of higher geometric structures, with applications to higher gauge theory \cite{Baez}.



There is a rich Lie theory relating Lie groupoids and algebroids, including subtleties not present in the usual
theory of Lie algebras and Lie groups. Every Lie groupoid $G\toto M$ gives  rise to a Lie algebroid $A_G\then M$, and  Lie's first theorem ({\em Lie I}) holds in this context: If a Lie algebroid $A$ admits an integration $G$, then there is a unique source-simply-connected one $G_A$ such that any other is the quotient $G_A\to G$ by a subgroupoid with discrete isotropies. The analogue of Lie's second theorem ({\em Lie II}) also holds, namely, if $G$ is source-simply-connected then every morphism of Lie algebroids $A_G\to A_H$ lifts to a groupoid map $G\to H$ in a unique way \cite{Mac-Xu2}. But Lie's third theorem ({\em Lie III}) does not hold in this generality: not every Lie algebroid is integrable. The precise obstructions to integrability were described in \cite{CF} in terms of the discreteness of certain monodromy groups.


The differentiation and integration procedures relating Lie groupoids and algebroids have several important applications;
in particular, understanding the global-infinitesimal correspondence of structures such as tensors, cohomologies, representations, actions and foliations, have drawn much recent attention, see e.g. \cite{AC0,AS2,bc,bcwz,CSS,ILX,jo,Mac-Xu,Mac-Xu2,mm1}. This includes far-reaching generalizations of the results relating Lie bialgebras and Poisson Lie groups, Poisson manifolds and symplectic groupoids, as well as the Van Est theorem in classical Lie theory.

In this brief note we outline the Lie theory of vector bundles over Lie groupoids and algebroids developed in \cite{BCdH} and discuss further applications to the study of more elaborate double structures \cite{mac-double1,mac-double2}.

\section{Vector bundles over Lie groupoids and Lie algebroids}


VB-groupoids and VB-algebroids \cite{Mac-book,GM} can be thought of as (categorified) vector bundles in the realm of Lie groupoids and Lie algebroids. They are simple examples of double structures, i.e., manifolds coupled with two different structures which are suitably compatible; these naturally arise e.g. in  Poisson Lie group theory \cite{mk2,mac-crelle}
(as Lie bialgebroids and double symplectic groupoids), or in the study of higher structures (see e.g. \cite{Baez} and references therein).

\begin{definition}\label{def:VB}
VB-groupoids and VB-algebroids are defined by diagrams as below, where $\toto$, $\then$ and $\to$ are used to denote Lie groupoid, Lie algebroid and vector bundle structures, respectively:
$$
\begin{matrix}
\Gamma & \toto & E\\
\downarrow & & \downarrow\\
G& \toto & M,
\end{matrix}
\qquad
\begin{matrix}
\Omega & \then & E\\
\downarrow & & \downarrow\\
A& \then & M.
\end{matrix}
$$
Here the vector bundle structural maps (projection, zero section, sum, scalar multiplication) are required to define Lie groupoid / Lie algebroid morphisms.
\end{definition}

Details on the definitions and other viewpoints can be found in \cite{GM,GM2}. One of our main contributions in \cite{BCdH} was a drastic simplification of the existing definitions, which allowed us to avoid dealing with several compatibility conditions when addressing the differentiation and integration of these objects.


VB-groupoids and VB-algebroids have been extensively studied in recent years, largely due to their ties with Poisson geometry and higher representation theory.
We mention here two fundamental examples of VB-groupoids -- there are analogous constructions for VB-algebroids.

\begin{example}
\
\begin{itemize}
\item The tangent and the cotangent bundles of a Lie groupoid $G\toto M$:
{\small{
$$
\begin{matrix}
TG & \toto & TM\\
\downarrow & & \downarrow\\
G& \toto & M,
\end{matrix}\qquad\qquad
\begin{matrix}
T^*G & \toto & A_G^*\\
\downarrow & & \downarrow\\
G& \toto & M.
\end{matrix}
$$}}
\item For a linear representation of a Lie groupoid $(G\toto M)$ on a vector bundle $(E\to M)$, the corresponding action groupoid is a VB-groupoid:
$$
\begin{matrix}
G\ltimes E & \toto & E\\
\downarrow & & \downarrow\\
G & \toto & M.
\end{matrix}
$$
Two extreme cases are when the groupoid is a manifold, and the representation is just a vector bundle, and when the groupoid is a group, in which case the notion of representation agrees with the usual one.
\end{itemize}
\end{example}


There is an important hidden piece of structure on VB-objects, the so-called {\em core}, which is a vector bundle $C\to M$ over the base manifold defined as follows: In the case of a VB-groupoid, the core is the kernel of the source map of $\Gamma$ restricted to the units $M$, and for VB-algebroids it is defined as the intersection of the two vector bundle projections defined on $\Omega$. The core of the tangent groupoid is $A_G\to M$ and the core of the cotangent groupoid is $T^*M$. There is an nontrivial duality theory for VB-objects, in which the roles of the core and side bundles are interchanged.

In the example of linear representations the core is 0 and, conversely, every VB-object with core 0 can be seen as a linear representation. In general, when nontrivial cores are present, VB-objects codify ``higher representations'': as proven in \cite{GM,GM2}, there is a one-to-one correspondence between isomorphism classes of VB-objects and of 2-term {\em representations up to homotopy}  \cite{AC1,AC2}; in the case of groupoids, these are certain non-associative actions, whereas in the case of algebroids they are described by non-flat connections. Representations up to homotopy are fundamental in the theory of Lie groupoids and algebroids, as they are needed to make sense of such basic objects as the adjoint representation.

In \cite{BCdH}, we establish the following fundamental infinitesimal-global correspondence for VB-objects;
its proof rests on a simplified approach to VB-objects, recalled in the next section.

\begin{theorem}\label{thm:main}
\
\begin{itemize}
\item {\bf Differentiation:}
Let $\Gamma\toto E$ be a VB-groupoid over $G\toto M$. Then the Lie algebroid $A_\Gamma\then E$ is naturally a VB-algebroid over $A_G\then M$.
\item {\bf Integration:}
Let $\Omega\then E$ be a VB-algebroid over $A\then M$.
If $\Omega$ is integrable, then its source-simply-connected integration $G_{\Omega}\toto M$ becomes a VB-groupoid over the source-simply connected integration $G_A\toto M$.
\end{itemize}
\end{theorem}


This result underlies the proofs of various theorems concerning the Lie theory of multiplicative structures on Lie groupoids, such as foliations \cite{jo} and more general Dirac structures \cite{ortiz}. It also leads to new results regarding the integration theory of representation up to homotopy, parallel to the formal viewpoint of \cite{AS2}.
The description of obstructions to integrability (in the sense of \cite{CF}) in the context of VB-algebroids is presented in \cite{BrCO}.


\section{Smooth vector bundles revisited}


In a smooth vector bundle, the operation of fibrewise addition turns out to be determined by scalar multiplication,
$x\mapsto h_\lambda(x)=\lambda x$; this follows from the fact that the map
$x\mapsto \frac{d}{d\lambda}\Big|_{\lambda=0}h_\lambda(x)$ identifies points of the total manifold with tangent vectors sitting at the zero section.
Of course not every action of the multiplicative monoid  $\Rm$ on a manifold $D$,
$h: \Rm \action D$, $h_1=Id$, $h_{\lambda\mu}=h_\lambda h_\mu$, corresponds to scalar multiplication relative to a vector bundle structure on $D$. A simple counter-example is given by $(\lambda,x)\mapsto \lambda^3x$ over the real line.
The following result provides a characterization of vector bundles by means of $(\R,\cdot)$-actions.

\begin{theorem}[\cite{GR}]\label{thm:GR}
A smooth action $(\R,\cdot)\action E$ is the scalar multiplication relative to a vector bundle structure on $E$
over the fixed points $h_0(E)$ if and only if it satisfies the regularity condition
$\frac{d}{d\lambda}\big |_{\lambda=0}h_\lambda(x)=0 \implies x=h_0(x)$.
\end{theorem}


This result allows us to rephrase the theory of vector bundles solely in terms of regular actions, which is convenient from a categorical viewpoint. For example, a vector bundle map is the same as an equivariant map, a subbundle is an invariant submanifold, etc. More importantly,  it leads to a simplified approach to more elaborate objects, such as double vector bundles -- which are described just as a pair of regular commuting actions \cite{GR}.


The previous characterization of vector bundles is an echo of a deeper phenomenon. The category of differential graded (super)manifolds can be seen to be equivalent to the category of manifolds equipped with an action of $\underline{\rm End}(\R^{0|1})$, the inner hom of endomorphisms of the odd line. The monoid $\Rm$ represents the body of $\underline{\rm End}(\R^{0|1})$ and is responsible for the grading -- and vector bundles are just graded manifolds generated in degrees 0 and 1. A first extension of the previous theorem along these lines is developed in \cite{BGG}.


In order to extend the above characterization of vector bundles to
the context of Lie groupoids and algebroids,  it is instructive to outline the main facts underpinning it:
out of any smooth action $h: \Rm\action D$, we obtain
\begin{itemize}
\item an embedded submanifold of fixed points, $M=h_0(D)\hookrightarrow D$,
\item a submersion $h_0: D \to M$,
\item a {\em vertical vector bundle} $V_h D: = \ker(d h_0)|_{M}
$, and
\item a {\em vertical lift} map $\mathcal{V}_h: D\to V_hD$,
$\mathcal{V}_h(x)=\frac{d}{d\lambda}\Big|_{\lambda=0}h_\lambda(x)$,
which is $\Rm$-equivariant.
\end{itemize}
\noindent
The key observation, leading to the proof of the previous theorem, is that the vertical lift map is a diffeomorphism if and only if the action is regular; in this case, it identifies $D$ with $V_hD$.


The construction of the fixed-point submanifold and of the vertical bundle can be naturally expressed as fibred-product
constructions. Thus, through a careful analysis of fibred products in the categories of Lie groupoids and algebroids,
one verifies that, when $D$ is replaced by a Lie groupoid or algebroid, the main facts above lead to the following simplified formulations of VB-groupoids and VB-algebroids \cite{BCdH}:

\begin{theorem}\label{thm:charac}\

\begin{itemize}
\item A VB-groupoid is the same as  a regular action $\Rm \action (\Gamma\toto E)$ by
Lie groupoid morphisms.
\item A VB-algebroid is the same as a regular action $\Rm \action (\Omega\then E)$ by
Lie algebroid morphisms.
\end{itemize}
\end{theorem}


Hence VB-groupoids and VB-algebroids are simply Lie groupoids and Lie algebroids equipped
with an additional compatible action of the monoid $\Rm$. From this viewpoint, the Lie-theoretic correspondence between VB-groupoids and VB-algebroids becomes a lot more transparent, since it boils down to the study of differentiation and integration of a single map (the action), that is covered by Lie II, along with an analysis of
the behavior of regularity of actions under these procedures.


\section{Double structures and the algebroid-Poisson duality}


Our results, described in Theorems \ref{thm:main} and \ref{thm:charac}, find applications in
the study of more general double structures -- those where the additional structure on a Lie groupoid or algebroid
is not just that of a vector bundle, but a Lie algebroid. The key to our approach is a duality that is central in Poisson geometry and extends the
correspondence between Lie algebra structures on a vector space $V$ and linear Poisson brackets on $C^\infty(V^*)$, which we now recall.


A Lie algebroid structure on a vector bundle $A\to M$ is the same as a Poisson structure on $A^*$ that is linear, in the sense that the bracket with a (fiberwise) linear function preserves linear and basic functions.
This dual viewpoint to Lie algebroids, in terms of Poisson structures, is often advantageous; e.g., Lie algebroid maps
(which can be complicated to work with) may be seen as linear maps whose dual relations are coisotropic submanifolds.
In terms of $\Rm$-actions $h$, linear Poisson structures are characterized by the fact that $(h_\lambda)_*\pi = \lambda \pi$, so we have the following equivalent approaches to Lie algebroids:
$$\begin{matrix} A \\ \Downarrow \\ M
\end{matrix}
\qquad \leftrightsquigarrow \qquad
\begin{matrix}(A^*,\pi)\\
\downarrow \\ M\end{matrix} \qquad \leftrightsquigarrow \qquad
\begin{matrix}\Rm\stackrel{h}{\action} (A^*,\pi)\\
{\scriptstyle{\pi\stackrel{h_\lambda}{\mapsto} \lambda \pi}}
\end{matrix}$$


In this spirit, by dualizing the double vector bundle underlying a VB-algebroid, one obtains a profitable
alternative characterization of VB-algebroids as objects dual to {\em double linear
Poisson structures}:
$$
\begin{matrix}
(\Omega^*_E,\pi) & \to & E\\
\downarrow & & \downarrow\\
C^*& \to & M.
\end{matrix}
$$
This viewpoint, as well as its formulation in terms of $\Rm$-actions, are important in establishing the algebroid counterpart of Theorem~\ref{thm:charac}.


Moving further into the world of double structures, we pass from VB-algebroids to ``LA-algebroids'', better known as \textit{double Lie algebroids} \cite{mac-double1,mk2,mac-double2}:
these are compatible diagrams of Lie algebroids,
$$
\begin{matrix}
\Omega & \then & E\\
\Downarrow & & \Downarrow\\
A& \then & M.
\end{matrix}
$$
Paradigmatic examples are the tangent bundle of a Lie algebroid and the cotangent bundle of a Lie bialgebroid.


Double Lie algebroids, in their original formulation \cite{mac-double1,mac-double2}, are rather intricate objects, involving subtle compatibility conditions. By resorting to the algebroid-Poisson duality, we obtained in \cite[Sec.~5.2]{BCdH} an easier-to-handle approach: double Lie algebroids are dual to {\em PVB-algebroids}, which are  VB-algebroids as in Def. \ref{def:VB} for which $\Omega$ is equipped with a Poisson structure $\pi$ which is linear with respect to $\Omega\to A$ and makes $(\Omega\then E, \pi)$ into a Lie bialgebroid. Making use of our previous characterization of VB-algebroids, one can go a step further and view PVB-algebroids as Lie bialgebroids $(\Omega\then E,\pi)$
endowed with a regular action $h: \Rm \action (\Omega\then E)$ such that $(h_\lambda)_*\pi = \lambda \pi$.


Double Lie algebroids arise as infinitesimal versions of double Lie groupoids (certain strict models of 2-groupoids). The problem of integrating a double Lie algebroid to a double Lie groupoid, however, is not fully understood yet (see \cite{luca}).
But this integration problem can be split into two steps, involving an intermediate object:
$$\text{Double Lie algebroids} \quad\leadsto\quad
\text{LA-groupoids} \quad\leadsto\quad
\text{Double Lie groupoids}$$
As an application of our results, we solved the first integration in this sequence \cite{BCdH}:

\begin{theorem}
If the total algebroid of a double Lie algebroid is integrable, then its source-simply-connected integration is an LA-groupoid.
\end{theorem}


An {\em LA-groupoid} \cite{mac-double1} is a compatible diagram of Lie groupoids and Lie algebroids:
$$\begin{matrix}
\Gamma & \toto & E\\
\Downarrow & & \Downarrow\\
G& \toto & M.
\end{matrix}$$
These objects play a role in modeling Lie algebroids over stacks \cite{waldron}. By means of the algebroid-Poisson duality, an LA-groupoid admits a simpler description as a VB-groupoid whose dual is equipped with a linear Poisson structure that is multiplicative (a Poisson groupoid); this dual object is referred to as a {\em PVB-groupoid} \cite{mk2}.


In \cite{BCdH}, we introduce PVB-algebroids and characterize both PVB-group\-oids and PVB-algebroids by using regular $\Rm$-actions; we prove that differentiation and integration behave naturally with respect to duality, and conclude that the Lie theory for LA-groupoids and double Lie algebroids, as in the previous theorem, is equivalent to the Lie theory of regular actions on Poisson groupoids and Lie bialgebroids,
$$
\Rm \action (\Gamma\toto E,\pi) \;\;\;
\stackrel{Lie}{\leftrightsquigarrow}\;\;\; \Rm \action (\Omega\then
E,\pi'),
$$ 
which we establish in \cite{BCdH}. The next integration, from LA- to double groupoids, involves additional topological issues yet to be clarified.






\begin{thebibliography}{99}

\bibitem{AC0}
Arias-Abad, C., Crainic, M.: The Weil algebra and the Van Est isomorphism.  {\em Ann. Inst. Fourier} {\bf 60} (2011), 927--970.

\bibitem{AC1}
Arias-Abad, C., Crainic, M.: Representations up to homotopy of Lie
algebroids.  {\em J. Reine Angew. Math.} {\bf 663} (2012), 91--126.

\bibitem{AC2}
Arias-Abad, C., Crainic, M.: Representations up to homotopy and
Bott's spectral sequence for Lie groupoids.  {\em Adv. Math.} {\bf 248}
(2013), 416--452.


\bibitem{AS2}
Arias-Abad, C., Sch\"atz, F.: The $A_\infty$ de Rham theorem and
integration of representations up to homotopy. {\em Int. Math. Res. Not.}
IMRN 2013, no. 16, 3790--3855.


\bibitem{Baez} Baez, J., Schreiber, U.: Higher gauge theory. Categories in Algebra,
Geometry and Mathematical Physics, eds. A. Davydov et al, Contemp. Math. 431, AMS, Providence, 2007.


\bibitem{BatesWein}
Bates, S., Weinstein, A.: {\em Lectures on the geometry of quantization. }
Berkeley Mathematics Lecture Notes, vol. 8.
American Mathematical Society, Providence, RI; Berkeley Center for
Pure and Applied Mathematics, Berkeley, CA, 1997.

\bibitem{BKS}
Bojowald, M., Kotov, A., Strobl, T.: Lie algebroid morphisms, Poisson
sigma models, and off-shell closed gauge symmetries. {\em J. Geom. Phys.} {\bf 54} (2005), 400--426.


\bibitem{BrCO}
Brahic, O., Cabrera, A., Ortiz, C.: Obstructions to the
integrability of VB-algebroids. Preprint arXiv:1403.1990 [math.DG].

\bibitem{BGG}
Bruce, A., Grabowska, K., Grabowski, J.: Graded bundles in the category of Lie groupoids.
{\em SIGMA} {\bf 11} (2015), 090, 25 pages.

\bibitem{bc}
Bursztyn, H.,  Cabrera, A.: Multiplicative forms at the
infinitesimal level. {\em Math. Ann.} {\bf 353} (2012), 663--705.


\bibitem{BCdH}
Bursztyn, H.,  Cabrera, A., del Hoyo, M.: Vector bundles over Lie groupoids and algebroids.
{\em Advances in Math.}, {\bf 290} (2016), 163--207.


\bibitem{bcwz}
Bursztyn, H., Crainic, M., Weinstein, A., Zhu, C., Integration of
twisted Dirac brackets, {\em Duke Math. J.} {\bf 123} (2004), 549-607.

\bibitem{CW}
Cannas da Silva, A., Weinstein, A., {\em Geometric models for
noncommutative algebras}. Berkeley Mathematics Lecture Notes, 10.
American Mathematical Society, Providence, RI; Berkeley Center for
Pure and Applied Mathematics, Berkeley, CA, 1999.

\bibitem{catfel}
Cattaneo, A., Felder, G., {\em Poisson sigma models and symplectic
groupoids}.  Quantization of singular symplectic quotients,  61--93,
Progr. Math., {\bf 198}, Birkhauser, Basel, 2001.

\bibitem{Connes}
Connes, A., {\em Noncommutative geometry}. Academic Press, San Diego, CA, 1994.

\bibitem{CDW}
Coste, A., Dazord, P., Weinstein, A., \textit{Groupo\"ides
symplectiques}. Publications du D\'epartement de Math\'ematiques.
Nouvelle S\'erie. A, Vol. 2,  i--ii, 1--62, Publ. D\'ep. Math.
Nouvelle S\'er. A, 87-2, Univ. Claude-Bernard, Lyon, 1987.


\bibitem{CF}
Crainic, M., Fernandes, R. L.: Integrability of Lie brackets. {\em
Annals of Math.} {\bf 157}  (2003), 575--620.

\bibitem{CSS}
Crainic, M., Salazar, M., Struchiner, I.: Multiplicative forms and Spencer operators.
{\em Math. Z.} {\bf 279} (2015), 939--979.




\bibitem{GR}
Grabowski, J., Rotkiewicz, M.: Higher vector bundles and
multi-graded symplectic manifolds. {\em J. Geom. Phys.} {\bf 59}
(2009), 1285--1305.



\bibitem{GM}
Gracia-Saz, A., Mehta, R.: Lie algebroid structures on double vector
bundles and representation theory of Lie algebroids. {\em Adv.
Math.} {\bf 223} (2010), 1236--1275.

\bibitem{GM2}
Gracia-Saz, A., Mehta, R.: $\mathcal{VB}$-groupoids and
representation theory of Lie groupoids. Preprint arxiv:1007.3658v4
[math.DG].

\bibitem{ILX}
Iglesias Ponte, D., Laurent-Gangoux, C., Xu, P.: Universal lifting
theorem and quasi-Poisson groupoids. {\em JEMS} {\bf 14} (2012), 681--731.

\bibitem{Hawkings}
Hawkings, E.: A groupoid approach to quantization. {\em J. Symplectic Geom.} {\bf 6} (2008), 61--125.


\bibitem{jo}
Jotz, M., Ortiz, C.: Foliated groupoids and infinitesimal ideal
systems. {\it Indag. Math.} {\bf 25} (2014), 1019--1053.

%



\bibitem{mac-double1}
Mackenzie, K.: Double Lie algebroids and second order geometry, I,
{\em Adv. Math.} {\bf 94} (1992), 180--239.


\bibitem{mk2}
Mackenzie, K.: On symplectic double groupoids and the duality of
Poisson groupoids. {\em Internat. J. Math.} {\bf 10} (1999),
435--456.

\bibitem{mac-double2}
Mackenzie, K.: Double Lie Algebroids and Second-Order Geometry, II,
{\em Adv. Math.} {\bf 154} (2000), 46--75.



\bibitem{Mac-book}
Mackenzie, K.: General theory of Lie groupoids and Lie algebroids;
London Math. Society Lecture Note Series 213, Cambridge University
Press, Cambridge, 2005.

\bibitem{mac-crelle}
Mackenzie, K.: Ehresmann doubles and Drinfel'd doubles for Lie
algebroids and Lie bialgebroids. {\em J. Reine Angew. Math.} {\bf
658} (2011), 193--245.

\bibitem{Mac-Xu}
Mackenzie, K.,  Xu, P.: Lie bialgebroids and Poisson groupoids. {\em
Duke Math. J.} {\bf 73} (1994), 415--452.

\bibitem{Mac-Xu2}
Mackenzie, K., Xu, P.: Integration of Lie bialgebroids.  {\em
Topology}  {\bf 39}  (2000), 445--467.




\bibitem{mm1}
Moerdijk, I., Mrcun, J.:
On integrability of infinitesimal actions.
{\em American J. Math.}
{\bf 124} (2002), 567--593.



\bibitem{ortiz}
Ortiz, C.: Multiplicative Dirac structures. {\em Pacific J. Math.}
{\bf 266} (2013), 329--365.

%



\bibitem{luca}
Stefanini, L.: On the integration of {LA}-groupoids and duality for
Poisson groupoids.{\em Travaux math\'ematiques}. {\bf 17} (2007),
65--85.



\bibitem{waldron}
Waldron, J.: Lie algebroids over differentiable stacks. {\em PhD thesis, University of York}, 2014. ArXiv:1511.07366.


\bibitem{we}
Weinstein, A.: Lagrangian mechanics and Lie groupoids. Mechanics day (Waterloo, ON, 1992), 207–231, {\em Fields Inst. Commun.}, 7, Amer. Math. Soc., Providence, RI, 1996.


\end{thebibliography}
\end{document}